\documentclass[a4paper,12pt]{report}
\usepackage[T2A]{fontenc}
\usepackage[cp1251]{inputenc}
\usepackage{amsfonts,amsmath,amssymb}
\usepackage[ukrainian,english,russian,english]{babel}
\usepackage{amsthm}
\usepackage{amscd}
\usepackage{amsmath}
\usepackage{amssymb}
\usepackage{graphicx}
\usepackage[T2A]{fontenc}
\usepackage[cp1251]{inputenc}
\def\beq{\begin{eqnarray}}
\def\eeq{\end{eqnarray}}

\theoremstyle{definition}

\theoremstyle{remark}

\begin{document}


{\bf Bounded solutions of dynamical systems in Hilbert space. }

{\bf А.А.Pokutnyi$^{*}$ \\
}

{\it * Institute of mathematics of NAS of Ukraine


e-mail: lenasas@gmail.com \\ }

Necessary and sufficient conditions for existence of bounded on the
entire real axis solutions of Schrodinger equation are obtained
under assumption that the homogeneous equation admits an exponential
dichotomy on the semi-axes. Bounded analytical solutions are
represented using generalized Green's operator. \\

{\bf Linear case.}

{\bf Statement of the problem.} Consider the next differential Shrodinger equation \cite{dal} \beq \label{11}
\frac{d\varphi(t)}{dt} = -i H(t)\varphi(t) + f(t), t \in J \eeq in a Hilbert space $H$, where, for each
 $t \in J \subset \mathbb{R}$, the unbounded operator $H(t)$ has the form $H(t) = H_{0} + V(t)$, where $H_{0} = H_{0}^{*}$ unbounded self-adjoint operator  with domain $D = D(H_{0}) \subset H;$ the mapping $t \rightarrow V(t)$ - is strongly continuous. Define as in $\cite{Rid2}$ operator-valued function
 $$
 \tilde{V}(t) = e^{itH_{0}}V(t)e^{-itH_{0}}.
 $$
 In this case for $\tilde{V}(t)$ Dyson's $\cite[p.311]{Rid2}$ representation is true and its propagator we define $\tilde{U}(t,s)$. If $U(t,s) = e^{-itH_{0}}\tilde{U}(t,s)e^{isH_{0}}$ then $\psi_{s}(t) = U(t,s)\psi$ is a weak solution of (\ref{11}) with condition $\varphi_{s}(s) = \psi$ in the sense that for any $\eta \in D(H_{0})$ function $(\eta, \psi_{s}(t))$ is differentiable and
 $$
 \frac{d}{dt}(\eta, \psi_{s}(t)) = -i(H_{0}\eta, \psi_{s}(t)) - i (V(t)\eta, \psi_{s}(t)), t \in J.
 $$

The present part dealt with the derivation of necessary and sufficient conditions for the existence of weak (in different senses) bounded and periodic solutions of the inhomogeneous equation (\ref{11}) with $f \in BC(J,H)  = \{ f : J \rightarrow H;$ the function $f$ is continuous and bounded $\}.$ Here the boundedness is treated in the sense that $|||f||| = sup_{t \in J}||f(t)|| < \infty$. For simplicity we suppose that $D$ dense in $H$. The operator $U(t,s)$ is a bounded linear operator for fixed $t,s$, and
since the set $D$ is dense in $H$, we find that it can be extended to the entire
space $H$ by continuity, which is assumed in forthcoming considerations. The extension of the family
of evolution operators to the entire space is denoted in the same way.

{ \bf 1. Bounded solutions.}
Throughout the following, we use the notion of exponential dichotomy in the sense of \cite[p.245 of the Russian translation]{dal}. It is of special interest to analyze the exponential dichotomy on the half-lines $\mathbb{R}_{s}^{-} = (-\infty, s]$ and $\mathbb{R}_{s}^{+} = [s; \infty)$. [In this case, the projection-valued functions defined on half-lines will be denoted by $P_{-}(t)$ for all $t \geq s$ and $P_{+}(t)$ for all $t < s$ with constants $M_{1}, \alpha_{1}$ and $M_{2}, \alpha_{2},$ respectively ($\alpha_{1}, \alpha_{2}$ - entropy or Lyapunov coefficients on the half-lines).] Most of the results obtained below follows directly from \cite{Pok}. The main result of this section can be stated as follows.

{\bf Theorem 1.} {\it Let $\{U(t,s), t \geq s \in \mathbb{R}\}$ be the family of strongly continuous evolution operators associated with equation (\ref{11}). Suppose that the following conditions are satisfied.

1. The operator $U(t,s)$ admits exponential dichotomy on the half-lines $\mathbb{R}_{0}^{+}$ and $\mathbb{R}_{0}^{-}$ with projection-valued operator-functions $P_{+}(t)$ and $P_{-}(t)$, respectively.

2. The operator $D = P_{+}(0) - (I - P_{-}(0))$ is generalized-invertible.

Then the following assertions hold.

1. There exist weak solutions of equation (\ref{11}) bounded on the entire line if and only if the vector function $f \in BC(\mathbb{R}, H)$ satisfies the condition
\beq \label{02}
\int_{-\infty}^{+\infty}H(t)f(t)dt = 0, \eeq
where $H(t) = \mathcal{P}_{N(D^{*})}P_{-}(0)U(0,t)$.

2. Under condition (\ref{02}), the weak solutions of (\ref{11}) bounded on the entire line have the form
\beq \label{03} x_{0}(t,c) = U(t,0)P_{+}(0)\mathcal{P}_{N(D)}c + (G[f])(t,0) \forall c\in H, \eeq
where
$$
(G[f])(t,s) = \left\{
\begin{array}{rcl}
\int_{s}^{t}T(t,\tau)P_{+}(\tau)f(\tau)d\tau -
\int_{t}^{+\infty}T(t,\tau)(I - P_{+}(\tau))f(\tau)d\tau + \\
+ T(t,s)P_{+}(s)D^{+}[\int_{s}^\infty
T(s,\tau)(I-P_{+}(\tau))f(\tau)d\tau +\\
+ \int_{-\infty}^{s} T(s,\tau)P_{-}(\tau)f(\tau)d\tau],
\hspace{0.2cm} t \geq s\\
\int_{-\infty}^{t} T(t,\tau)P_{-}(\tau)f(\tau)d\tau - \int_t^s
T(t,\tau)(I-P_{-}(\tau))f(\tau)d\tau + \\
+ T(t,s)(I-P_{-}(s))D^{+}[\int_{s}^\infty
T(s,\tau)(I-P_{+}(\tau))f(\tau)d\tau +\\
+ \int_{-\infty}^{s} T(s,\tau)P_{-}(\tau)f(\tau)d\tau],
\hspace{0.2cm} s \geq t
\end{array}
\right.
$$
is the generalized Green operator of the problem on the bounded, on the entire line, solutions $$
(G[f])(0+,0) - (G[f])(0-,0) = -\int_{-\infty}^{+\infty}H(t)f(t)dt;
$$
$$
\mathcal{L}(G[f])(t,0) = f(t), ~~t\in \mathbb{R}
$$
and
$$
(\mathcal{L}x)(t) = \frac{dx}{dt} - iH(t)x(t),
$$
$D^{+}$ is the Moore-Penrouse pseudoinverse operator to the operator $D$; $\mathcal{P}_{N(D)} = I - D^{+}D$ and $\mathcal{P}_{N(D^{*})} = I - DD^{+}$ are the projections \cite{BoiSam} onto the kernel and cokernel of the operator $D$.}

{\bf Remark 1.} A similar theorem holds for the case in which the family of evolution operators $U(t,s)$ admits exponential dichotomy on the half-lines $\mathbb{R}_{s}^{+}$ and $\mathbb{R}_{s}^{-}$.

The proof of this theorem is the same as in \cite[Theorem 1]{Pok}.

Now we show that condition 2 in theorem 1 can be omitted and in the different senses equation (\ref{11}) is always resolvable without condition (\ref{02}).
From the proof of the theorem 1 follows that equation (\ref{11}) have bounded solutions if and only if the operator equation
\beq \label{03}
D\xi = g, \eeq
$$g = \int_{-\infty}^{0}U(0,\tau)P_{-}(\tau)f(\tau)d\tau + \int_{0}^{+\infty}U(0,\tau)(I - P_{+}(\tau))f(\tau)d\tau
$$
is resolvable and its number depends from the dimension of $N(D)$.

Consider next 3 cases.

1) Classical strong generalized solutions.

Consider case when the set of values of $D$ is closed ($R(D) = \overline{R(D)}$).
Then \cite{BoiSam} $g \in R(D)$ if and only if $\mathcal{P}_{N(D^{*})}g = 0$ and the set of solutions of (\ref{03}) has the form \cite{BoiSam} $\xi = D^{+}g + \mathcal{P}_{N(D)}c, \forall c \in H $.

2) Strong generalized solutions.
Consider the case when $R(D) \neq \overline{R(D)}$. We show that operator $D$ may be extended to $\tilde{D}$ in such way that $R(\tilde{D})$ is closed.

Since the operator $D$ is bounded the next representations  $H$ in the direct sum are true
$$
H = N(D) \oplus X, H = \overline{R(D)} \oplus Y,
$$
with $X = N(D)^{\bot}$ and $Y = \overline{R(D)}^{\bot}$.
Let $E = H/N(D)$ is quotient space of $H$ and $\mathcal{P}_{\overline{R(D)}}$ - orthoprojector, which projects onto $\overline{R(D)}$. Then operator
$$
D\mathcal{P}_{\overline{R(D)}} : E \rightarrow R(D) \subset \overline{R(D)},
$$
is linear, continuous and injective. In this case \cite[p.26,29]{Lash} we can define strong generalized solution of equation
$$
D\mathcal{P}_{\overline{R(D)}} \tilde{\xi} = g, x \in E .
$$
Fill up the space $E$ in the norm $||x||_{\overline{E}} = ||D\mathcal{P}_{\overline{R}(D)}x||_{F},$ where $F = \overline{R(D)}$ \cite{Lash}. Then extended operator
$$\overline{D\mathcal{P}_{\overline{R(D)}}} : \overline{E} \rightarrow \overline{R(D)}, E \subset \overline{E}$$ is homeomorphism of $\overline{E}$ and $\overline{R(D)}$.
By virtue of construction of strong generalized solution \cite{Lash} equation
$$
\overline{D\mathcal{P}_{\overline{R(D)}}} \overline{\xi} = g,
$$
has a unique generalized solution, which we denote $\tilde{D}^{+}g$.

3) Pseudosolutions.

Consider element $g \notin \overline{R(D)}$. This condition is equivalent $\mathcal{P}_{N(D^{*})}g \neq 0$. In this case there exists elements from $H$ which minimise norm $||D\xi - g ||_{H}$ :
$$
\xi = D^{+}g + \mathcal{P}_{N(D)}c, \forall c \in H.
$$
These elements are called pseudosolutions \cite{BoiSam}.

{\bf Remark 2.} It should be noted that in every case 1) - 3)  the form of bounded solutions (\ref{03}) isn't change.

As follows from 1) and 3) the notion of exponential dichotomy is
equivalent of existence of bounded on the entire real axis
solutions of (\ref{11}).

{\bf Main result (Nonlinear case).}

In the Hilbert space $H$, consider the differential equation
\beq \label{04}
\frac{d\varphi(t)}{dt} = -iH(t)\varphi(t) + \varepsilon Z(\varphi,t,\varepsilon) + f(t).
\eeq
We seek a bounded solution $\varphi(t,\varepsilon)$ of equation (\ref{04}) that becomes one of the solutions of the generating equation (\ref{11}) for $\varepsilon = 0$.

To find a necessary condition on the operator function $Z(\varphi, t, \varepsilon),$ we impose the joint constraints
$$
Z(\cdot, \cdot, \cdot) \in BC(\mathbb{R},H)\times C[0,\varepsilon_{0}]\times C[||x - x_{0}||\leq q],
$$
where $q$ is some positive constant.

Let us show that this problem can be solved with the use of the operator equation for generating constants
\beq \label{05} F(c) = \int_{-\infty}^{+\infty}H(t)Z(\varphi_{0}(t,c),t,0)dt = 0. \eeq

{\bf Theorem 2} {\it (necessary condition). Let the equation (\ref{11}) admit exponential dichotomy on the half-lines $\mathbb{R}_{0}^{+}$ and $\mathbb{R}_{0}^{-}$ with projection-valued operator functions $P_{+}(t)$ and $P_{-}(t)$, respectively, and let the nonlinear equation (\ref{04}) have a bounded solution $\varphi(\cdot, \varepsilon)$ that becomes one of the solutions of the generating equation (\ref{11}) with constant $c = c^{0},$ $\varphi(t, 0) = \varphi_{0}(t, c^{0})$ for $\varepsilon = 0$. Then this constant should satisfy the equation for generating constants (\ref{05}). }

The proof of this theorem is the same as in \cite[Theorem 1]{Pok}.

To find a sufficient condition for the existence of bounded solutions of (\ref{11}), we additionally assume that the operator function $Z(\varphi, t, \varepsilon)$ is strongly differentiable in a neighborhood of the generating solution $(Z(\cdot, t, \varepsilon) \in C^{1}[||x - x_{0}|| \leq q])$.

This problem can be solved with the use of the operator
$$
B_{0} = \int_{-\infty}^{+\infty}H(t)A_{1}(t)T(t,0)P_{+}(0)\mathcal{P}_{N(D)}dt : H \rightarrow H,
$$
where $A_{1}(t) = Z^{1}(v, t, \varepsilon)|_{v = \varphi_{0}, \varepsilon = 0}$ (the Fr\'{e}chet derivative).

{\bf Theorem 3} {\it (sufficient condition). Suppose that the equation (\ref{11}) admits exponential dichotomy on the half-lines $\mathbb{R}_{0}^{+}$ and $\mathbb{R}_{0}^{-}$ with projection-valued functions
$P_{+}(t)$ and $P_{-}(t),$ respectively. In addition, let the operator $B_{0}$ satisfy the following conditions.

1. The operator $B_{0}$ is pseudoinvertible.

2. $\mathcal{P}_{N(B_{0}^{*})}\mathcal{P}_{N(D^{*})}P_{-}(0) = 0$.

Then for an arbitrary element of $c = c^{0} \in H$ satisfying the equation for generating constants (\ref{05}). This solution can be found with the use of the iterative process

$$
\overline{y}_{k+1}(t, \varepsilon) = \varepsilon G[Z(\varphi_{0}(\tau, c^{0} + y_{k}, \tau, \varepsilon))](t,0),
$$
$$
c_{k} = -B_{0}^{+}\int_{-\infty}^{+\infty}H(\tau)\{A_{1}(\tau)\overline{y}_{k}(\tau, \varepsilon) + \mathcal{R}(y_{k}(\tau, \varepsilon), \tau, \varepsilon) \}d\tau,
$$
$$
y_{k+1}(t, \varepsilon) = T(t,0)P_{+}(0)\mathcal{P}_{N(D)}c_{k} + \overline{y}_{k+1}(t,0,\varepsilon),
$$
$$
x_{k}(t,\varepsilon) = x_{0}(t, c^{0}) + y_{k}(t,\varepsilon),~~~ k = 0,1,2,...,~~~y_{0}(t,\varepsilon) = 0, ~~~x(t,\varepsilon) = \lim_{k \rightarrow \infty} x_{k}(t, \varepsilon).
$$
}
{\bf Relationship between necessary and sufficient conditions.}

First, we formulate the following assertion.

{\bf Corollary.} {\it Let a functional $F(c)$ have the Fr\'{e}chet derivative $F^{(1)}(c)$ for each element $c^{0}$ of the Hilbert space $H$ satisfying the equation for generating constants (\ref{05}). If $F^{1}(c)$ has a bounded inverse, then equation (\ref{04}) has a unique bounded solution on the entire line for each $c^{0}$.}

{\bf Remark 3.} If assumptions of the corollary are satisfied, then it follows from its proof that the operators $B_{0}$ and $F^{(1)}(c^{0})$ are equal. Since the operator $F^{(1)}(c)$ is invertible, it follows that assumptions 1 and 2 of Theorem 3 are necessarily satisfied for the operator $B_{0}$. In this case, equation (\ref{04}) has a unique bounded solution for each $c^{0} \in H$. Therefore, the invertibility condition for the operator $F^{1}(c)$ relates the necessary and sufficient conditions. In the finite-dimensional case, the condition of invertibility of the operator $F^{(1)}(c)$ is equivalent to the condition of simplicity of the root $c^{0}$ of the equation for generating amplitudes \cite{BoiSam}.

 In such way we generalize the well-known method of Lyapunov-Schmidt. It should be emphasized that theorem 2 and 3 give us condition of chaotic behavior of (\ref{04}) \cite{Chuesh}.

\end{document}